\documentclass{amsart}
\usepackage{latexsym,amsmath,amssymb}
\theoremstyle{definition}
\theoremstyle{remark}

\numberwithin{equation}{section}
\begin{document}

\title[Weighted conditional type operators ]
{ Some Classes of weighted conditional type operators and their
spectra }

\author{\sc\bf Y. Estaremi  }
\address{\sc Y. Estaremi  }
\email{yestaremi@pnu.ac.ir}

\address{Department of mathematics, Payame Noor university , p. o. box: 19395-3697, Tehran,
Iran}

\thanks{}

\thanks{}



\date{}

\dedicatory{}

\commby{}

%%% ----------------------------------------------------------------------
\begin{abstract}
In this paper, some $\ast$-classes of weighted conditional
expectation type operators, such as $A$-class, $\ast$-$A$-class
and quasi-$\ast$-$A$ classes on $L^{2}(\Sigma)$ are investigated.
Also, the spectrum, point spectrum and spectral radius of these
operators are computed.

\noindent {}
\end{abstract}

\maketitle

\section{ \sc\bf Introduction and Preliminaries}
This paper is about an important operator in statistics and
analysis, that is called conditional expectation. By the
projection theorem, the conditional expectation E(X) is the best
mean square predictor of X in range E, also in analysis it is
proved that lots of operators are of the form E and of the form of
combinations of E and multiplications operators. I think our
informatioms in this paper and refereces [4,5] will be useful for
experts in statistics.
 Conditional expectations have been studied in an operator
theoretic setting, by, for example, de Pagter and Grobler [8], de
Pagter et al. [5], and Rao [18,19], as positive operators acting
on Lp-spaces or Banach function spaces. In [15] we showed that
many of the fundamental results in stochastic processes can be
formulated and proved in the measure-free framework of Riesz
spaces.

Let $(X,\Sigma,\mu)$ be a complete $\sigma$-finite measure space.
For any sub-$\sigma$-finite algebra $\mathcal{A}\subseteq
 \Sigma$, the $L^2$-space
$L^2(X,\mathcal{A},\mu_{\mid_{\mathcal{A}}})$ is abbreviated  by
$L^2(\mathcal{A})$, and its norm is denoted by $\|.\|_2$. All
comparisons between two functions or two sets are to be
interpreted as holding up to a $\mu$-null set. The support of a
measurable function $f$ is defined as $S(f)=\{x\in X; f(x)\neq
0\}$. We denote the vector space of all equivalence classes of
almost everywhere finite valued measurable functions on $X$ by
$L^0(\Sigma)$.

\vspace*{0.3cm} For a sub-$\sigma$-finite algebra
$\mathcal{A}\subseteq\Sigma$, the conditional expectation operator
associated with $\mathcal{A}$ is the mapping $f\rightarrow
E^{\mathcal{A}}f$, defined for all non-negative measurable
function $f$ as well as for all $f\in L^2(\Sigma)$, where
$E^{\mathcal{A}}f$, by the Radon-Nikodym theorem, is the unique
$\mathcal{A}$-measurable function satisfying
$$\int_{A}fd\mu=\int_{A}E^{\mathcal{A}}fd\mu, \ \ \ \forall A\in \mathcal{A} .$$
As an operator on $L^{2}({\Sigma})$, $E^{\mathcal{A}}$ is
idempotent and $E^{\mathcal{A}}(L^2(\Sigma))=L^2(\mathcal{A})$.
This operator will play a major role in our work. Let $f\in
L^0(\Sigma)$, then $f$ is said to be conditionable with respect to
$E$ if $f\in\mathcal{D}(E):=\{g\in L^0(\Sigma): E(|g|)\in
L^0(\mathcal{A})\}$. Throughout this paper we take $u$ and $w$ in
$\mathcal{D}(E)$. If there is no possibility of confusion, we
write $E(f)$ in place of $E^{\mathcal{A}}(f)$.  A detailed
discussion about this operator may be found in \cite{rao}.

\vspace*{0.2cm} Composition of conditional expectation operators
and multiplication operators  appear often in the study of other
operators such as multiplication operators and weighted
composition operators. Specifically, in \cite{mo}, S.-T. C. Moy
characterized all operators on $L^p$ of the form $f\rightarrow
E(fg)$ for $g$ in $L^q$ with $E(|g|)$ bounded. Eleven years later,
R. G. Douglas, \cite{dou}, analyzed positive projections on
$L^{1}$ and many of his characterizations are in terms of
combinations of multiplications and conditional expectations. More
recently, P.G. Dodds, C.B. Huijsmans and B. De Pagter, \cite{dhd},
extended these characterizations to the setting of function ideals
and vector lattices. J. Herron presented some assertions about the
operator $EM_u$ on $L^p$ spaces in \cite{her}. Also, some results
about multiplication conditional expectation operators can be
found in \cite{g, lam}. In \cite{e,ej} we investigated some
classic properties of multiplication conditional expectation
operators $M_wEM_u$ on $L^p$ spaces.\\

Let $\mathcal{H}$ be the infinite dimensional complex Hilbert
space and let $\mathcal{L(H)}$ be the algebra of all bounded
operators on $\mathcal{H}$. An operator $T\in \mathcal{L(H)}$ is a
partial isometry if $\|Th\|=\|h\|$ for $h$ orthogonal to the
kernel of $T$. It is known that an operator $T$ on a Hilbert space
is partial isometry if and only if $TT^{\ast}T=T$. Every operator
$T$ on a Hilbert space $\mathcal{H}$ can be decomposed into $T =
U|T|$ with a partial isometry $U$, where $|T| =
(T^*T)^{\frac{1}{2}}$ . $U$ is determined uniquely by the kernel
condition $\mathcal{N}(U) = \mathcal{N}(|T|)$. Then this
decomposition is called the polar decomposition. The Aluthge
transformation $\widehat{T}$ of the operator $T$ is defined by
$\widehat{T}=|T|^{\frac{1}{2}}U|T|^{\frac{1}{2}}$. The operator
$T$ is said to be positive operator and write $T\geq 0$, if
$\langle
Th, h\rangle\geq 0$, for all $h\in \mathcal{H}$. \\
In this paper we will be concerned with characterizing weighted
conditional expectation type operators on $L^2(\Sigma)$ in terms
of membership in the partial $A$-classes and computing the
spectrum, point spectrum, spectral radius and we prove that the
point spectrum and joint point spectrum of the quasi-$\ast$-$A$-
classes weighted conditional type operators are the same. Here is
a brief review of what constitutes membership for an operator $T$
on
a Hilbert space in each classes:\\

(i) $T$ belongs to $\ast$-$A$-class, if
$|T^2|\geq|T^{\ast}|^2$.\\

(ii) $T$ belongs to quasi-$\ast$-$A$-class, if
$T^{\ast}|T^2|T\geq T^{\ast}|T^{\ast}|^2T$.\\

(iii) The operator $T$ belongs to $A$-class, if $|T|^2\leq|T^2|$.\\

\section{ \sc\bf Some classes of weighted conditional type operators}

In the first we reminisce some theorems that we have proved in
\cite{ej}.

\vspace*{0.3cm} {\bf Theorem 2.1.}  The operator $T=M_wEM_u$ is
bounded on $L^{2}(\Sigma)$ if and only if
$(E|w|^{2})^{\frac{1}{2}}(E|u|^{2})^{\frac{1}{2}} \in
L^{\infty}(\mathcal{A})$, and in this case its norm is given by
$\|T\|=\|(E(|w|^{2}))^{\frac{1}{2}}(E(|u|^{2}))^{\frac{1}{2}}\|_{\infty}$.\\

\vspace*{0.3cm} {\bf Lemma 2.2.} Let $T=M_wEM_u$ be a bounded
operator on $L^{2}(\Sigma)$ and let $p\in (0,\infty)$. Then
$$(T^{\ast}T)^{p}=M_{\bar{u}(E(|u|^{2}))^{p-1}\chi_{S}(E(|w|^{2}))^{p}}EM_{u}$$
and
$$(TT^{\ast})^{p}=M_{w(E(|w|^{2}))^{p-1}\chi_{G}(E(|u|^{2}))^{p}}EM_{\bar{w}},$$
where $S=S(E(|u|^2))$ and $G=S(E(|w|^2))$.

\vspace*{0.3cm} {\bf Theorem 2.3.} The unique polar decomposition
of bounded operator $T=M_wEM_u$ is $U|T|$, where

$$|T|(f)=\left(\frac{E(|w|^{2})}{E(|u|^{2})}\right)^{\frac{1}{2}}\chi_{S}\bar{u}E(uf)$$
and
 $$U(f)=\left(\frac{\chi_{S\cap
 G}}{E(|w|^{2})E(|u|^{2})}\right)^{\frac{1}{2}}wE(uf),$$
for all $f\in L^{2}(\Sigma)$.

\vspace*{0.3cm} {\bf Theorem 2.4.} The Aluthge transformation of
$T=M_wEM_u$ is
$$\widehat{T}(f)=\frac{\chi_{S}E(uw)}{E(|u|^{2})}\bar{u}E(uf), \ \ \ \ \ \  \ \ \ \  \ \  \ f\in L^{2}(\Sigma).$$

\vspace*{0.3cm} {\bf Remark 2.5.} By Theorem 2.3 and Theorem 2.4,
we can compute the polar decomposition and Aluthge transformation
of $T^\ast=U^\ast|T^\ast|$ as follows:

\begin{align*}
|T^\ast|(f)&=\left(\frac{E(|u|^{2})}{E(|w|^{2})}\right)^{\frac{1}{2}}\chi_{G}\bar{w}E(\bar{w}f);\\
U^\ast(f)&=\left(\frac{\chi_{S\cap G}}{(E(|u|^{2})E(|w|^{2})}\right)^{\frac{1}{2}}\bar{u}E(\bar{w}f);\\
\widehat{T^\ast}(f)&=\frac{\chi_{G}E(\overline{uw})}{E(|w|^{2})}wE(\bar{w}f),
\end{align*}
for all $f\in L^{2}(\Sigma)$.\\

In the sequel some necessary and sufficient conditions for
weighted conditional type operators $M_wEM_u$ to be in
$\ast$-$A$-classes, quasi-$\ast$-$A$-classes and $A$-classes of operators, will be presented.\\

\vspace*{0.3cm} {\bf Theorem 2.6.}
 Let $T=M_wEM_u$
be a bounded operator on $L^{2}(\Sigma)$. Then\\

(a) If $|E(uw)|^2\geq E(|u|^2)E(|w|^2)$ on $S=S(E(|u|^2))$, then
$T$ belongs to $A$-class.\\

(b) If $T$ belongs to $A$-class, then $|E(uw)|^2\geq
E(|u|^2)E(|w|^2)$ on $S'=S(E(u))$.\\

(c) If $S=S'$, then $T$ belongs to $A$-class if and only if
$|E(uw)|^2\geq E(|u|^2)E(|w|^2)$ on $S$.

\vspace*{0.3cm} {\bf Proof.} (a) By Theorem 2.3 we have

$$|T|^2(f)=E(|w|^2)\chi_{S}\bar{u}E(uf), \ \ \ \
|T^2|(f)=|E(uw)|(\frac{E(|w|^2)}{E(|u|^2)})^{\frac{1}{2}}\chi_{S}\bar{u}E(uf),$$
for all $f\in L^2(\Sigma)$. So for every $f\in L^2(\Sigma)$ we
have

$$\langle|T^2|(f)-|T|^2(f),f\rangle=\int_{X}|E(uw)|(\frac{E(|w|^2)}{E(|u|^2)})^{\frac{1}{2}}\chi_{S}\bar{uf}E(uf)-E(|w|^2)\chi_{S}\bar{uf}E(uf)d\mu$$
$$=\int_{X}|E(uw)|(\frac{E(|w|^2)}{E(|u|^2)})^{\frac{1}{2}}\chi_{S}|E(uf)|^2-E(|w|^2)\chi_{S}|E(uf)|^2d\mu.$$
This implies that if $|E(uw)|^2\geq E(|u|^2)E(|w|^2)$ on $S$, then
$\langle|T^2|(f)-|T|^2(f),f\rangle\geq0$ for all $f\in
L^2(\Sigma)$. So $T$ is belonged to $A$-class.\\

(b) If $T$ belongs to $A$-class, then by (a), for all $f\in
L^2(\Sigma)$ we have

$$\int_{X}|E(uw)|(\frac{E(|w|^2)}{E(|u|^2)})^{\frac{1}{2}}\chi_{S}|E(uf)|^2-E(|w|^2)\chi_{S}|E(uf)|^2d\mu\geq0.$$

Let $A\in \mathcal{A}$, with $0<\mu(A)<\infty$. By replacing $f$
to $\chi_{A}$, we have

$$\int_{A}|E(uw)|(\frac{E(|w|^2)}{E(|u|^2)})^{\frac{1}{2}}\chi_{S}|E(u)|^2-E(|w|^2)\chi_{S}|E(u)|^2d\mu\geq0.$$

Since $A\in \mathcal{A}$ is arbitrary, then $|E(uw)|^2\geq
E(|u|^2)E(|w|^2)$ on $S'=S(E(u))$.\\

(c) It follows from (a) and (b).$\Box$

\vspace*{0.3cm} {\bf Theorem 2.7.} Let $T=M_wEM_u$ be a bounded
operator on $L^{2}(\Sigma)$. Then\\

(a) If
$$u|E(uw)|^{\frac{1}{2}}(\frac{E(|w|^2)}{E(|u|^2)})^{\frac{1}{4}}\chi_{S}\geq
\bar{w}(E(|u|^2))^{\frac{1}{2}},$$ then $T$ belongs to
$\ast$-$A$-class.\\

(b) If $T$ is belonged to $\ast$-$A$-class, then
$$|E(u)|^2|E(uw)|(\frac{E(|w|^2)}{E(|u|^2)})^{\frac{1}{2}}\chi_{S}\geq
(E(|u|^2))^{\frac{1}{2}}|E(w)|^2.$$

{\bf Proof.} (a) By using Theorem 2.3,  for every $f\in
L^2(\Sigma)$

$$|T^2|(f)=|E(uw)|(\frac{E(|w|^2)}{E(|u|^2)})^{\frac{1}{2}}\chi_{S}\bar{u}E(uf), \ \ \ \ |T^{\ast}|^2(f)=E(|u|^2)wE(\bar{w}f),$$
hence for every $f\in L^2(\Sigma)$

$$\langle|T^2|(f)-|T^{\ast}|^2(f),f\rangle=\int_{X}|E(uw)|(\frac{E(|w|^2)}{E(|u|^2)})^{\frac{1}{2}}\chi_{S}\bar{uf}E(uf)-E(|u|^2)w\bar{f}E(\bar{w}f)d\mu$$
$$\int_{X}|E(uw)|(\frac{E(|w|^2)}{E(|u|^2)})^{\frac{1}{2}}\chi_{S}|E(uf)|^2-E(|u|^2)|E(\bar{w}f)|^2d\mu.$$

This implies that, if
$$u|E(uw)|^{\frac{1}{2}}(\frac{E(|w|^2)}{E(|u|^2)})^{\frac{1}{4}}\chi_{S}\geq
\bar{w}(E(|u|^2))^{\frac{1}{2}},$$ then $|T^2|\geq|T^{\ast}|^2$
i.e, $T$ belongs to $\ast$-$A$-class.\\

(b) If $T$ belongs to $\ast$-$A$-class, then by (a), for all $f\in
L^2(\Sigma)$ we have

$$\int_{X}|E(uw)|(\frac{E(|w|^2)}{E(|u|^2)})^{\frac{1}{2}}\chi_{S}|E(uf)|^2-E(|u|^2)|E(\bar{w}f)|^2d\mu\geq0.$$

Let $A\in \mathcal{A}$, with $0<\mu(A)<\infty$. By replacing $f$
to $\chi_{A}$, we have

$$\int_{A}|E(uw)|(\frac{E(|w|^2)}{E(|u|^2)})^{\frac{1}{2}}\chi_{S}|E(uf)|^2-E(|u|^2)|E(\bar{w}f)|^2d\mu\geq0.$$

Since $A\in \mathcal{A}$ is arbitrary, then

$$|E(u)|^2|E(uw)|(\frac{E(|w|^2)}{E(|u|^2)})^{\frac{1}{2}}\chi_{S}\geq
(E(|u|^2))^{\frac{1}{2}}|E(w)|^2.$$

\vspace*{0.3cm} {\bf Theorem 2.8.}
 Let $T=M_wEM_u$
be a bounded operator on $L^{2}(\Sigma)$. Then\\

(a) If $|E(uw)|^2\geq E(|u|^2)E(|w|^2)$, then $T$ belongs to
quasi-$\ast$-$A$-class.\\

(b) If $T$ belongs to quasi-$\ast$-$A$-class, then
$$|E(uw)|^3(E(|w|^2))^{\frac{1}{2}}\geq
(E(|u|^2))^{\frac{3}{2}}(E(|w|^2))^2.$$

 \vspace*{0.3cm} {\bf Proof.} (a) Direct computation shows that, for
$f\in L^2(\Sigma)$,

$$T^{\ast}|T^2|T(f)=|E(uw)|^3(\frac{E(|w|^2)}{E(|u|^2)})^{\frac{1}{2}}\chi_{S}\bar{u}E(uf), \ \ \ T^{\ast}|T^{\ast}|^2T(f)=E(|u|^2)(E(|w|^2))^2\bar{u}E(uf).$$
So, for all $f\in L^2(\Sigma)$

$$\langle T^{\ast}|T^2|T(f)-T^{\ast}|T^{\ast}|^2T(f),f\rangle$$$$
=\int_{X}|E(uw)|^3(\frac{E(|w|^2)}{E(|u|^2)})^{\frac{1}{2}}\chi_{S}\bar{uf}E(uf)-E(|u|^2)(E(|w|^2))^2\bar{uf}E(uf)d\mu$$
$$=\int_{X}(|E(uw)|^3(\frac{E(|w|^2)}{E(|u|^2)})^{\frac{1}{2}}\chi_{S}-E(|u|^2)(E(|w|^2))^2)|E(uf)|^2d\mu.$$
This implies that if $|E(uw)|^2\geq E(|u|^2)E(|w|^2)$, then
$$T^{\ast}|T^2|T-T^{\ast}|T^{\ast}|^2T\geq0.$$\\

(b) If $T$ is belonged to quasi-$\ast$-$A$-class, then by (a), for
all $f\in L^2(\Sigma)$ we have

$$\int_{X}|E(uw)|(\frac{E(|w|^2)}{E(|u|^2)})^{\frac{1}{2}}\chi_{S}|E(uf)|^2-E(|u|^2)|E(\bar{w}f)|^2d\mu\geq0.$$

Let $A\in \mathcal{A}$, with $0<\mu(A)<\infty$. By replacing $f$
to $\chi_{A}$, we have

$$\int_{A}|E(uw)|(\frac{E(|w|^2)}{E(|u|^2)})^{\frac{1}{2}}\chi_{S}|E(u)|^2-E(|u|^2)|E(\bar{w})|^2d\mu\geq0.$$

Since $A\in \mathcal{A}$ is arbitrary, then

$$|E(uw)|^3(E(|w|^2))^{\frac{1}{2}}\geq
(E(|u|^2))^{\frac{3}{2}}(E(|w|^2))^2.$$\\

J. Herron showed that:
 $EM_u$ is normal if and only if $u\in
L^{\infty}(\mathcal{A})$ \cite{her}. So we have the following
result.

\vspace*{0.3cm} {\bf Corollary 2.9.} If we consider the operator $T=EM_u:L^2(\Sigma)\rightarrow L^2(\Sigma)$, then the followings are equivalent.\\

(a) $T$ is normal.\\

(b) $T$ is belonged to quasi-$\ast$-$A$-class.\\

(c) $u\in L^{\infty}(\mathcal{A})$.\\

\vspace*{0.3cm} {\bf Example 2.10.} {\bf (a)} Let $X=[0,1]\times
[0,1]$, $d\mu=dxdy$, $\Sigma$  the  Lebesgue subsets of $X$ and
let $\mathcal{A}=\{A\times [0,1]: A \ \mbox{is a Lebesgue set in}
\ [0,1]\}$. Then, for each $f$ in $L^2(\Sigma)$, $(Ef)(x,
y)=\int_0^1f(x,t)dt$, which is independent of the second
coordinate. This example is due to A. Lambert and B. Weinstock
\cite{la}. Now, if we take $u(x,y)=y^{\frac{x}{8}}$ and $w(x,
y)=\sqrt{(4+x)y}$, then $E(|u|^2)(x,y)=\frac{4}{4+x}$ and
$E(|w|^2)(x,y)=\frac{4+x}{2}$. So, $E(|u|^2)(x,y)E(|w|^2)(x,y)=2$
and $|E(uw)|^2(x,y)=64\frac{4+x}{(x+12)^2}$. Direct computations
shows that $E(|u|^2)(x,y)E(|w|^2)(x,y)\leq|E(uw)|^2(x,y)$. Thus,
by Theorem 2.6 the weighted conditional type operator $M_wEM_u$
belongs to $A$-classes of operators on $L^2(\Sigma)$.

 {\bf (b)} Let $\Omega=[-1,1]$, $d\mu=\frac{1}{2}dx$ and $\mathcal{A}=<\{(-a,a):0\leq a\leq1\}>$ (Sigma algebra generated by symmetric intervals).
Then
 $$E^{\mathcal{A}}(f)(x)=\frac{f(x)+f(-x)}{2}, \ \ x\in \Omega,$$
 where $E^{\mathcal{A}}(f)$ is defined. Let $u(x)=x^2-1$ and $w\equiv1$, then by Theorem 2.8 the operator
 $M_wEM_u$ belongs to quasi-$\ast$-$A$-classes of operators on
 $L^2(\Sigma)$.\\

 \section{ \sc\bf The spectra of weighted conditional type operators}

In this section, we shall denote by $\sigma(T)$, $\sigma_{p}(T)$,
$\sigma_{jp}(T)$, $r(T)$ the spectrum of $T$, the point spectrum
of $T$, the joint point spectrum of $T$, the spectral radius of
$T$, respectively. The spectrum of an operator $T$ is the set
$$\sigma(T)=\{\lambda\in \mathbb{C}:T-\lambda I \ \ \  {\text{is \  not \
invertible}}\}.$$ A complex number $\lambda\in \mathbb{C}$ is said
to be in the point spectrum $\sigma_{p}(T)$ of the operator $T$,
if there is a unit vector $x$ satisfying $(T-\lambda)x=0$. If in
addition, $(T^{\ast}-\bar{\lambda})x=0$, then $\lambda$ is said to
be in the joint spectrum $\sigma_{jp}(T)$ of $T$. Also, the spectral radius of $T$ is defined by $r(T)=\sup\{|\lambda|: \lambda\in \sigma(T)\}$.\\
If $A, B\in \mathcal{B}(\mathcal{H})$, then it is well known that
$\sigma(AB)\setminus\{0\}=\sigma(BA)\setminus\{0\}$ and
$\sigma_{p}(AB)\setminus\{0\}=\sigma_{p}(BA)\setminus\{0\}$.
 J. Herron
showed that if $EM_u:L^2(\Sigma)\rightarrow L^2(\Sigma)$, then
$\sigma(EM_u)=ess \ range(E(u))\cup\{0\}$, \cite{her}.\\

 Let
$A_{\lambda}=\{x\in X:E(u)(x)=\lambda\}$, for $\lambda\in
\mathbb{C}$. Suppose that $\mu(A_{\lambda})>0$. Since
$\mathcal{A}$ is $\sigma$-finite,  there exists an
$\mathcal{A}$-measurable subset $B$ of $A_{\lambda}$ such that
$0<\mu(B)<\infty$, and $f=\chi_{B}\in L^p(\mathcal{A})\subseteq
L^p(\Sigma)$. Now
$$EM_u(f)-\lambda f=E(u)\chi_{B}-\lambda \chi_{B}=0.$$ This
implies that $\lambda\in \sigma_{p}(EM_u)$.\\
If there exists $f\in L^p(\Sigma)$ such that $f\chi_{C}\neq 0$
$\mu$-a.e, for $C\in \Sigma$ of positive measure and
$E(uf)=\lambda f$ for $0\neq \lambda \in \mathbb{C}$, then
$f=\frac{E(uf)}{\lambda}$, which means that $f$ is
$\mathcal{A}$-measurable. Therefore $E(uf)=E(u)f=\lambda f$ and
$(E(u)-\lambda)f=0$. This implies that $C\subseteq A_{\lambda}$
and so $\mu(A_{\lambda})>0$. We have proved the following theorem.

\vspace*{0.3cm} {\bf Theorem 3.1.} Let
$T=EM_u:L^2(\Sigma)\rightarrow L^2(\Sigma)$. Then\\

(a)
$$\sigma_{p}(EM_u)\setminus\{0\}=\{\lambda\in\mathbb{C}\setminus\{0\}:\mu(A_{\lambda})>0\}$$

(b)
$$ \{\lambda\in\mathbb{C}:\mu(A_{\lambda})>0\}\subseteq\sigma_{p}(EM_u).$$

(c) If $\mu(\{x\in X:E(u)(x)=0\})>0$, then
$$ \{\lambda\in\mathbb{C}:\mu(A_{\lambda})>0\}=\sigma_{p}(EM_u).$$

\vspace*{0.3cm} {\bf Theorem 3.2.} Let
$T=M_wEM_u:L^2(\Sigma)\rightarrow L^2(\Sigma)$. Then\\

(a) $$\sigma(M_wEM_u)\setminus \{0\}=ess \
range(E(uw))\setminus\{0\}.$$

(b) If $S\cap G=X$, then $$\sigma(M_wEM_u)=ess \ range(E(uw)),$$
where $S=S(E(|u|^2))$ and $G=S(E(|w|^2))$.

(c)
$$\sigma_{p}(M_wEM_u)\setminus\{0\}=\{\lambda\in\mathbb{C}\setminus\{0\}:\mu(A_{\lambda,w})>0\},$$
where $A_{\lambda,w}=\{x\in X:E(uw)(x)=\lambda\}$.\\

 \vspace*{0.3cm} {\bf Proof.} (a) Since
 $$\sigma(M_wEM_u)\setminus\{0\}=\sigma(EM_uM_w)\setminus\{0\}=\sigma(EM_{uw})\setminus\{0\}=ess
 \ range(E(uw))\setminus\{0\},$$
 then we have
 $$\sigma(M_wEM_u)\setminus \{0\}=ess \
range(E(uw))\setminus\{0\}.$$

(b) We know that $\sigma(EM_{uw})=ess \ range(E(uw))$. So, we have
to prove that $0\notin\sigma(EM_{uw})$ if and only if
$0\notin\sigma(M_wEM_u)$ by (a).\\
Let $0\notin\sigma(EM_{uw})$. Then $EM_{uw}$ is surjective and so
$\mathcal{A}=\Sigma$. Thus $E=I$. So $0\notin\sigma(M_wEM_u)$.\\
Conversely, by Theorem 2.3 the polar decomposition of
$M_wEM_u=U|M_wEM_u|$ is as follow,
$$|M_wEM_u|(f)=\left(\frac{E(|w|^{2})}{E(|u|^{2})}\right)^{\frac{1}{2}}\chi_{S}\bar{u}E(uf)$$
and
 $$U(f)=\left(\frac{\chi_{S\cap
 G}}{E(|w|^{2})E(|u|^{2})}\right)^{\frac{1}{2}}wE(uf),$$
for all $f\in L^{2}(\Sigma)$.\\

If $0\notin\sigma(M_wEM_u)$, then $|M_wEM_u|$ is invertible and
$U$ is unitary. Therefore $U^{\ast}U=UU^{\ast}=I$. The equation
$UU^{\ast}=I$ implies that $w\in L^{0}(\mathcal{A}_{S\cap G})$,
where $\mathcal{A}_{S\cap G}=\{A\cap S\cap G: \ A\in
\mathcal{A}\}$. Since $S\cap G=X$, then $w\in L^{0}(\mathcal{A})$.
Hence $0\notin\sigma(M_wEM_u)=\sigma(EM_{uw})$.\\

(c) By Theorem 3.1 we have
$$\sigma_{p}(M_wEM_u)\setminus\{0\}=\sigma_{p}(EM_{u}M_{w})\setminus\{0\}=\sigma_{p}(EM_{uw})\setminus\{0\}.$$
So
$$\sigma_{p}(M_wEM_u)\setminus\{0\}=\{\lambda\in\mathbb{C}\setminus\{0\}:\mu(A_{\lambda,w})>0\}.$$\\

For each natural number $n$, we define
$$\triangle_n(T)=\widehat{\triangle_{n-1}{T}} \ \ \ \ \
\triangle_1(T)=\triangle(T)=\widehat{T}.$$

We call $\triangle_n(T)$ the $n$-th Aluthge transformation of $T$.
It is proved that $r(T)=\lim_{n\rightarrow
\infty}\|\triangle_n(T)\|$ in \cite{ty}. So, we have the following theorems.\\

\vspace*{0.3cm} {\bf Theorem 3.3.} For
$T=M_wEM_u:L^2(\Sigma)\rightarrow L^2(\Sigma)$, we have
$r(T)=\|E(uw)\|_{\infty}$.\\

\vspace*{0.3cm} {\bf Proof.} By Theorem 2.4 we conclude that for
every natural number $n$ we have
$\triangle_n(T)=\triangle(T)=\widehat{T}$. Hence
$r(T)=\|\widehat{T}\|=\|E(uw)\|_{\infty}$.\\

%\vspace*{0.3cm} {\bf Theorem 3.4.} If $|E(uw)|^2\geq
 %E(|u|^2)E(|w|^2)$, then\\

% (a) If $T=M_wEM_u$ does not have a dense range, then
 %$T\mid_{\overline{ran T}}$ is belonged to $\ast$-$A$-class.\\

% (b) If $M$ is an invariant subspace of $T=M_wEM_u$, then
 %$T\mid_{M}$ is belonged to
% quasi-$\ast$-$A$-class.\\

\vspace*{0.3cm} {\bf Theorem 3.4.}
 If $T=M_wEM_u$ belongs to quasi-$\ast$-$A$-class, then
 $$\sigma_{p}(M_wEM_u)=\sigma_{jp}(M_wEM_u).$$

 \vspace*{0.3cm} {\bf Proof.} Let $f\in L^2(\Sigma)\setminus \{0\}$ and $\lambda\in \mathbb{C}$, such that
 $wE(uf)=\lambda f$. Let $M=span\{f\}$, the closed linear subspace
 generated by $f$. Thus we can represent $T=M_wEM_u$ as the
 following $2\times 2$ operator matrix with respect to the
 decomposition $L^2(\Sigma)=M\oplus M^{\perp}$,
 $$
 T= \left[
         \begin{array}{rr}
              M_{\lambda} & PM_{w}EM_{u}- M_{\lambda} \\
              0 &  M_{w}EM_{u}-PM_{w}EM_{u}
          \end{array} \right].
$$
where $P$ is is the orthogonal projection of $L^2(\Sigma)$ onto
$M$. Since $T$ belongs to quasi-$\ast$-$A$-class, for every $g\in
L^2(\Sigma)$, there exists $\alpha\in \mathbb{C}$ such that
$$P(g)=\alpha \lambda f=\alpha T(f)=T(\alpha f),$$
and
$$<P|T^2|P(g)-P|T^{\ast}|^2P(g),g>=<|T^2|P(g)-|T^{\ast}|^2P(g),Pg>$$$$=<T^{\ast}|T^2|T(\alpha f)-T^{\ast}|T^{\ast}|^2T(\alpha
f),\alpha f>\geq0.$$

Hence $P(|T^2|-|T^{\ast}|^2)P\geq0$. Direct computation shows that
$$P|T^2|^2P= \left[
         \begin{array}{rr}
              M_{|\lambda|^4} & 0 \\
              0 &  0
          \end{array} \right]
.$$

Thus
$$\left[
         \begin{array}{rr}
              M_{|\lambda|^2} & 0 \\
              0 &  0
          \end{array} \right]=(P|T^2|^2P)^{\frac{1}{2}}\geq
          P|T^2|P\geq P|T^{\ast}|^2P $$$$=PTT^{\ast}P=\left[
         \begin{array}{rr}
              M_{|\lambda|^2}+AA^{\ast} & 0 \\
              0 &  0
          \end{array} \right]
,$$ where $A=PM_{w}EM_{u}- M_{\lambda}$. This implies that $A=0$
and so $M$ is invariant under $T^{\ast}=M_{\bar{u}}EM_{\bar{w}}$,
too. Hence we conclude that
$T^{\ast}(f)=M_{\bar{u}}EM_{\bar{w}}(f)=\lambda f$. This means
that $\sigma_{p}(M_wEM_u)=\sigma_{jp}(M_wEM_u)$.\\

\vspace*{0.3cm} {\bf Corollary 3.5.} If $|E(uw)|^2\geq
 E(|u|^2)E(|w|^2)$, then\\

 (a)
 $$\sigma_{jp}(M_wEM_u)\setminus \{0\}=ess \
 range(E(uw))\setminus\{0\}=\{\lambda\in\mathbb{C}\setminus\{0\}:\mu(A_{\lambda,w})>0\}.$$

(b) If $S\cap G=X$, then

 $$\sigma_{jp}(M_wEM_u)=ess \
 range(E(uw))=\{\lambda\in\mathbb{C}:\mu(A_{\lambda,w})>0\}.$$

%If $S_1=S(u(E(|w|^2))^{\frac{1}{2}}$, then $L^2(X\backslash
%S_1)\subseteq N(T)$. Also, if
%$S_2=S(\bar(w)(E(|u|^2))^{\frac{1}{2}}$, then $L^2(X\backslash
%S_2)\subseteq N(T^{\ast})$. Thus, if $S_2\neq X$, then $T$ isn't
%dense range. Let $P$ be an orthogonal projection onto
%$$\overline{R(T)}$. Hence $PT=T$, $T^{\ast}P=T^{\ast}$, $PTP=TP$
%and $PT^{\ast}=PT^{\ast}P$.

%The operator $T$ is $M$-paranormal, if for all unit vectors $x$ in
%$H$, $\|Tx\|^2\leq M\|T^2x\|$. $T$ is paranormal if and only if
%$T^{\ast^2}T^2-2kT^{\ast}T+k^2\geq0$, $k\in \mathbb{R}$ and $T$ is
%$M$-paranormal if
%and only if $M^2T^{\ast^2}T^2+2kT^{\ast}T+k^2\geq0$.\\

%If $T$ is paranormal, then
%$$E(|w|^2)|E(u)|^2\leq |E(uw)|(E(|w|^2))^{\frac{1}{2}}|E(u)|.$$

\end{document}